\documentclass[oneside,english]{amsart}
\usepackage[T1]{fontenc}
\usepackage[latin9]{inputenc}
\usepackage{amstext}
\usepackage{amsthm}
\usepackage{amssymb}
\usepackage{babel}
\begin{document}
\title{Closed Form Equations for Triangular Numbers Multiple of Other Triangular
Numbers}
\author{Vladimir PLETSER}
\address{European Space Agency (ret.)}
\begin{abstract}
Triangular numbers that are multiple of other triangular numbers are
investigated. It is known that for any positive non-square integer
multiplier, there is an infinity of multiples of triangular numbers
which are triangular numbers. If the multiplier is a squared integer,
there is either one or no solution, depending on the multiplier value.
Instead of recurrent relations, we develop in this paper closed form
equations to calculate directly the values of triangular numbers and
their indices without the need of knowing the previous solutions.
We develop the theoretical equations for four cases of ranks from
1 to 4 and we give several examples for non-square multipliers 2,
3, 5 and 8.
\end{abstract}

\maketitle
\textbf{Keywords}: Triangular Numbers, Multiple of Triangular Numbers,
Closed Form Equations

\section{Introduction}

Triangular numbers are defined as $T_{t}=\frac{t\left(t+1\right)}{2}$
and enjoy many properties, relations and formulas (see e.g. \cite{key-1,key-2}).
Triangular numbers $T_{\xi}$ that are multiples of other triangular
number $T_{t}$ 
\begin{equation}
T_{\xi}=kT_{t}\label{eq:1}
\end{equation}
have been investigated in the past in some specific cases (\cite{key-3,key-4,key-5,key-6,key-7,key-8}).
Recently, Pletser showed (\cite{key-9}) that, first, for square integer
values of $k$, there are either no or only one solution of (\ref{eq:1}),
depending on the value of $k$; and second, for non-square integer
values of $k$, there are infinitely many solutions of (\ref{eq:1})
and recurrent relations can be found for the four variables $t,\xi,Tt$
and $T_{\xi}$

\begin{align}
t_{n} & =2\left(\kappa+1\right)t_{n-r}-t_{n-2r}+\kappa\label{eq:3.3}\\
\xi_{n} & =2\left(\kappa+1\right)\xi_{n-r}-\xi_{n-2r}+\kappa\label{eq:3.3-1}\\
T_{t_{n}} & =\left(4\left(\kappa+1\right)^{2}-2\right)T_{t_{n-r}}-T_{t_{n-2r}}+\left(T_{\kappa}-\gamma\right)\label{eq:3.3-2}\\
T_{\xi_{n}} & =\left(4\left(\kappa+1\right)^{2}-2\right)T_{\xi_{n-r}}-T_{\xi_{n-2r}}+k\left(T_{\kappa}-\gamma\right)\label{eq:3.3-3}
\end{align}
These recurrent relations involve three parameters specific to each
case of the multiplier $k$, i.e., a rank $r$, defined as the number
of successive values of $t$ solutions of (\ref{eq:1}) with slowly
decreasing ratios of the $r^{\text{th}}$ to the $\left(r-1\right)^{\text{th}}$
values of $t$, $t_{r-1}/t_{r}$. The two other parameters are $\kappa$
and $\gamma$ respectively the sum $\kappa=t_{r-1}+t_{r}$ and the
product $\gamma=t_{r-1}t_{r}$ of the first two sequential values
of $t_{n}$ for $n=r-1$ and $r$. Several relations exist between
these parameters and are investigated in \cite{key-9}. Note that
only cases with $k>1$ are of interest as $k=0$ and $k=1$ yield
obvious solutions respectively $\xi=0$ and $\xi=t$, both $\forall t$.
Although triangular numbers $T_{t}$ are usually defined for $t\in\mathbb{Z}^{+}$,
triangular numbers can be extended to negative indices $t<0$ as $T_{-t}=T_{t-1}$.

In this paper, instead of solutions with recurrent relations, we investigate
closed form solutions of (\ref{eq:1}) for the case of non-square
integer values of $k$.

\section{Closed Form Equations}

\subsection{General equations}

One wishes to calculate directly the $n^{\text{th}}$value of $t_{n},\xi_{n},T_{t_{n}}$
and $T_{\xi_{n}}$ without necessarily knowing smaller values needed
for recurrent solutions. Closed form equations are easy to calculate
from the recursive relations (\ref{eq:3.3}) to (\ref{eq:3.3-3}),
which are all linear and non-homogeneous. Let us consider first (\ref{eq:3.3}).
The associated closed form equation is the sum of the homogeneous
solution $t_{h_{n}}$ and a particular solution $t_{p}$ (for more
details, see e.g. \cite{key-10,key-11,key-12}).

The homogeneous characteristic equation associated to (\ref{eq:3.3})
reads successively
\begin{align}
x^{2r}-2\left(\kappa+1\right)x^{r}+1 & =0\label{eq:5-4}\\
\left(x^{r}-\alpha^{r}\right)\left(x^{r}-\beta^{r}\right) & =0\label{eq:5-5}
\end{align}
with
\begin{align}
\alpha & =\left(\left(\kappa+1\right)+\sqrt{\kappa\left(\kappa+2\right)}\right)^{\frac{1}{r}}\label{eq:5-6}\\
\beta & =\left(\left(\kappa+1\right)-\sqrt{\kappa\left(\kappa+2\right)}\right)^{\frac{1}{r}}\label{eq:5-7}
\end{align}
with the obvious relation $\left(\alpha\beta\right)^{r}=1$. This
equation (\ref{eq:5-5}) has $2r$ distinct characteristic roots,
including real and complex roots and, for odd $r$, $r=2\rho+1$,
one has
\begin{align}
x^{r}-\alpha^{r} & =\left(x-\alpha\right)\left(\sum_{j=0}^{2\rho}\alpha^{j}x^{2\rho-j}\right)\nonumber \\
 & =\left(x-\alpha\right)\prod_{j=1}^{\rho}\left[x-\alpha\left(\cos\left(\frac{2j\pi}{r}\right)\pm i\sin\left(\frac{2j\pi}{r}\right)\right)\right]\label{eq:5-8}
\end{align}
and for even $r$, $r=2\rho$,
\begin{align}
x^{r}-\alpha^{r} & =\left(x-\alpha\right)\left(x+\alpha\right)\left(\sum_{j=0}^{\rho-1}\alpha^{2j}x^{2\left(\rho-j-1\right)}\right)\nonumber \\
 & =\left(x-\alpha\right)\left(x+\alpha\right)\prod_{j=1}^{\rho-1}\left[x-\alpha\left(\cos\left(\frac{j\pi}{\rho}\right)\pm i\sin\left(\frac{j\pi}{\rho}\right)\right)\right]\label{eq:5-9}
\end{align}
and similar relations for $\beta$ replacing $\alpha$. Recall that
complex characteristic roots introduce the element of periodicity.
The homogeneous solution is therefore respectively (\ref{eq:5-10})
for odd $r$, $r=2\rho+1$ and (\ref{eq:5-11}) for even $r$, $r=2\rho$,
\begin{align}
t_{h_{n}} & =\left(A+\sum_{j=1}^{\rho}\left[A_{j}^{\prime}\cos\left(\frac{2nj\pi}{r}\right)+A_{j}^{\prime\prime}\sin\left(\frac{2nj\pi}{r}\right)\right]\right)\alpha^{n}\label{eq:5-10}\\
 & +\left(B+\sum_{j=1}^{\rho}\left[B_{j}^{\prime}\cos\left(\frac{2nj\pi}{r}\right)+B_{j}^{\prime\prime}\sin\left(\frac{2nj\pi}{r}\right)\right]\right)\beta^{n}\nonumber \\
t_{h_{n}} & =\left(A+\left(-1\right)^{n}A^{\prime}+\sum_{j=1}^{\rho-1}\left[A_{j}^{\prime\prime}\cos\left(\frac{nj\pi}{\rho}\right)+A_{j}^{\prime\prime\prime}\sin\left(\frac{nj\pi}{\rho}\right)\right]\right)\alpha^{n}\label{eq:5-11}\\
 & +\left(B+\left(-1\right)^{n}B^{\prime}+\sum_{j=1}^{\rho}\left[B_{j}^{\prime\prime}\cos\left(\frac{nj\pi}{\rho}\right)+B_{j}^{\prime\prime\prime}\sin\left(\frac{nj\pi}{\rho}\right)\right]\right)\beta^{n}\nonumber 
\end{align}
with the $2r$ constants $A,A^{\prime},A_{j}^{\prime\prime},A_{j}^{\prime\prime\prime},B,B^{\prime},B_{j}^{\prime\prime}$
and $B_{j}^{\prime\prime\prime}$ to be determined by $2r$ boundary
conditions, namely the first $2r$ values of $t_{n}$, for $n=0$
to $n=2r-1$. In practice, one can simplify this method by using the
fact that $t_{j}=t_{-\left(j+1\right)}$, i.e. for $n=-r$ to $n=r-1$,
necessitating to only know the first $\left(r-1\right)$ values of
$t_{n}$.

The particular solution $t_{p}$ can be found by posing it equal to
a constant $t_{p}=\tau$. This particular solution must be a solution
of the initial equation (\ref{eq:3.3}). Therefore replacing $t_{n},t_{n-r}$
and $t_{n-2r}$ by $t_{p}=\tau$ in (\ref{eq:3.3}) yields directly
$\tau=-1/2$. The complete solution is then $t_{n}=t_{h_{n}}-1/2$.

\subsection{Closed form equations for $r=1$ ($\rho=0$)}

The homogeneous solution (\ref{eq:5-10}) reduces to
\begin{equation}
t_{h_{n}}=A\alpha^{n}+B\beta^{n}
\end{equation}
yielding the complete solution 
\begin{equation}
t_{n}=A\alpha^{n}+B\beta^{n}-1/2
\end{equation}
with the boundary conditions $t_{0}=t_{-1}=0$, yielding $A+B=1/2$
and $A\alpha^{-1}+B\beta^{-1}=1/2$, giving $A=\alpha\left(1-\beta\right)/2\left(\alpha-\beta\right)$
and $B=\beta\left(\alpha-1\right)/2\left(\alpha-\beta\right)$.

\subsubsection{Closed forms for $k=2$}

The smallest case of $k$ having a rank equal to unity is $k=2$,
for which $t_{r-1}=0$, $t_{r}=2$ (see OEIS \cite{key-13}, A053141)
and $\kappa=2$, yielding from (\ref{eq:5-6}) and (\ref{eq:5-7}),
$\alpha=\left(3+2\sqrt{2}\right)=\left(1+\sqrt{2}\right)^{2}$, $\beta=\left(3-2\sqrt{2}\right)=\left(1-\sqrt{2}\right)^{2}$,
and further, $A=\left(2+\sqrt{2}\right)/8$, $B=\left(2-\sqrt{2}\right)/8$,
and the closed form equation
\begin{equation}
t_{n}=\frac{1}{8}\left(\left(2+\sqrt{2}\right)\left(1+\sqrt{2}\right)^{2n}+\left(2-\sqrt{2}\right)\left(1-\sqrt{2}\right)^{2n}\right)-\frac{1}{2}
\end{equation}
The closed form equations for the three other variables for this case
of $k=2$ are found similarly:

- for $\xi_{n}$ (\ref{eq:3.3-1}) (see OEIS \cite{key-13}, A001652)
\begin{equation}
\xi_{n}=\frac{1}{4}\left(\left(1+\sqrt{2}\right)^{2n+1}+\left(1-\sqrt{2}\right)^{2n+1}\right)-\frac{1}{2}
\end{equation}
- for $T_{t_{n}}$ (\ref{eq:3.3-2}) (see OEIS \cite{key-13}, A075528),
with $\left(17\pm12\sqrt{2}\right)=\left(1\pm\sqrt{2}\right)^{4}$,
\begin{equation}
T_{t_{n}}=\frac{1}{64}\left(\left(1+\sqrt{2}\right)^{4n+2}+\left(1-\sqrt{2}\right)^{4n+2}\right)-\frac{3}{32}
\end{equation}
- for $T_{\xi_{n}}$ (\ref{eq:3.3-3}) (see OEIS \cite{key-13}, A029549)
\begin{equation}
T_{\xi_{n}}=\frac{1}{32}\left(\left(1+\sqrt{2}\right)^{4n+2}+\left(1-\sqrt{2}\right)^{4n+2}\right)-\frac{3}{16}
\end{equation}

\subsubsection{Closed forms for $k=3$}

For $k=3$, $r=1$ and $t_{r-1}=0$, $t_{r}=1$ (see OEIS \cite{key-13},
A061278) and $\kappa=1$, yielding $\alpha=\left(2+\sqrt{3}\right)$,
$\beta=\left(2-\sqrt{3}\right)$ from (\ref{eq:5-6}) and (\ref{eq:5-7}),
and $A=\left(3+\sqrt{3}\right)/12$, $B=\left(3-\sqrt{3}\right)/12$,
and the closed form
\begin{equation}
t_{n}=\frac{1}{12}\left(\left(3+\sqrt{3}\right)\left(2+\sqrt{3}\right)^{n}+\left(3-\sqrt{3}\right)\left(2-\sqrt{3}\right)^{n}\right)-\frac{1}{2}
\end{equation}
Similarly, the closed forms for the three other variables are:

- for $\xi_{n}$ (\ref{eq:3.3-1}) (see OEIS \cite{key-13}, A001571):
\begin{equation}
\xi_{n}=\frac{1}{4}\left(\left(1+\sqrt{3}\right)\left(2+\sqrt{3}\right)^{n}+\left(1-\sqrt{3}\right)\left(2-\sqrt{3}\right)^{n}\right)-\frac{1}{2}
\end{equation}
- for $T_{t_{n}}$ (\ref{eq:3.3-2}) (see OEIS \cite{key-13}, A076139),
with $\left(7\pm4\sqrt{3}\right)=\left(2\pm\sqrt{3}\right)^{2}$:
\begin{equation}
T_{t_{n}}=\frac{1}{48}\left(\left(2+\sqrt{3}\right)^{2n+1}+\left(2-\sqrt{3}\right)^{2n+1}\right)-\frac{1}{12}
\end{equation}
- for $T_{\xi_{n}}$ (\ref{eq:3.3-3}) (see OEIS \cite{key-13}, A076140):
\begin{equation}
T_{\xi_{n}}=\frac{1}{16}\left(\left(2+\sqrt{3}\right)^{2n+1}+\left(2-\sqrt{3}\right)^{2n+1}\right)-\frac{1}{4}
\end{equation}

\subsection{Closed form equations for $r=2$ ($\rho=1$)}

The homogeneous solution (\ref{eq:5-11}) reduces to 
\begin{equation}
t_{h_{n}}=\left(A+\left(-1\right)^{n}A^{\prime}\right)\alpha^{n}+\left(B+\left(-1\right)^{n}B^{\prime}\right)\beta^{n}
\end{equation}
yielding the complete solution 
\begin{equation}
t_{n}=\left(A+\left(-1\right)^{n}A^{\prime}\right)\alpha^{n}+\left(B+\left(-1\right)^{n}B^{\prime}\right)\beta^{n}-1/2
\end{equation}
with the boundary conditions $t_{0}=t_{-1}=0$ and $t_{1}=t_{-2}$,
yielding the four relations 
\begin{eqnarray*}
A+A^{\prime}+B+B^{\prime} & = & 1/2\\
\left(A-A^{\prime}\right)\alpha^{-1}+\left(B-B^{\prime}\right)\beta^{-1} & = & 1/2\\
\left(A-A^{\prime}\right)\alpha+\left(B-B^{\prime}\right)\beta & = & t_{1}+1/2\\
\left(A+A^{\prime}\right)\alpha^{-2}+\left(B+B^{\prime}\right)\beta^{-2} & = & t_{1}+1/2
\end{eqnarray*}
yielding
\begin{align*}
A & =\frac{\alpha\left(\alpha+1\right)\left(1-\beta^{2}\right)+2t_{1}\left(\alpha-1\right)}{4\left(\alpha^{2}-\beta^{2}\right)}\\
A^{\prime} & =\frac{\alpha\left(\alpha-1\right)\left(1-\beta^{2}\right)-2t_{1}\left(\alpha+1\right)}{4\left(\alpha^{2}-\beta^{2}\right)}\\
B & =\frac{\beta\left(\beta+1\right)\left(\alpha^{2}-1\right)-2t_{1}\left(\beta-1\right)}{4\left(\alpha^{2}-\beta^{2}\right)}\\
B^{\prime} & =\frac{\beta\left(\beta-1\right)\left(\alpha^{2}-1\right)+2t_{1}\left(\beta+1\right)}{4\left(\alpha^{2}-\beta^{2}\right)}
\end{align*}

\subsubsection{Closed forms for $k=5$}

The smallest case of $k$ with rank two is $k=5$ (see OEIS \cite{key-13},
A077259), for which $\kappa=8$ yielding $\alpha=\sqrt{9+4\sqrt{5}}=2+\sqrt{5}$
and $\beta=\sqrt{9-4\sqrt{5}}=2-\sqrt{5}$, with $t_{1}=2$. The complete
solution is then

\begin{align}
t_{n} & =\frac{1}{20}\left(\left(\left(5+2\sqrt{5}\right)-\left(-1\right)^{n}\sqrt{5}\right)\left(2+\sqrt{5}\right)^{n}\right.\\
 & \left.+\left(\left(5-2\sqrt{5}\right)+\left(-1\right)^{n}\sqrt{5}\right)\left(2-\sqrt{5}\right)^{n}\right)-\frac{1}{2}\nonumber 
\end{align}
which yields for even $n$
\begin{equation}
t_{n}=\frac{1}{20}\left(\left(5+\sqrt{5}\right)\left(2+\sqrt{5}\right)^{n}+\left(5-\sqrt{5}\right)\left(2-\sqrt{5}\right)^{n}\right)-\frac{1}{2}
\end{equation}
and for odd $n$
\begin{equation}
t_{n}=\frac{1}{20}\left(\left(5+3\sqrt{5}\right)\left(2+\sqrt{5}\right)^{n}+\left(5-3\sqrt{5}\right)\left(2-\sqrt{5}\right)^{n}\right)-\frac{1}{2}
\end{equation}

\subsubsection{Closed forms for $k=8$}

For $k=8$ (see OEIS \cite{key-13}, A336623), $r=2$, for which $\kappa=16$
yielding $\alpha=\sqrt{17+12\sqrt{2}}=\left(1+\sqrt{2}\right)^{2}$
and $\beta=\sqrt{17-12\sqrt{2}}=\left(1-\sqrt{2}\right)^{2}$, with
$t_{1}=5$. The complete solution is then 
\begin{align}
t_{n} & =\frac{1}{8}\left(\left(\frac{3}{2}\left(2+\sqrt{2}\right)-\left(-1\right)^{n}\left(1+\sqrt{2}\right)\right)\left(1+\sqrt{2}\right)^{2n}\right.\\
 & \left.+\left(\frac{3}{2}\left(2-\sqrt{2}\right)-\left(-1\right)^{n}\left(1-\sqrt{2}\right)\right)\left(1-\sqrt{2}\right)^{2n}\right)-\frac{1}{2}\nonumber 
\end{align}
which can be dissociated for even $n$, yielding
\begin{equation}
t_{n}=\frac{1}{16}\left(\left(4+\sqrt{2}\right)\left(1+\sqrt{2}\right)^{2n}+\left(4-\sqrt{2}\right)\left(1-\sqrt{2}\right)^{2n}\right)-\frac{1}{2}
\end{equation}
and for odd $n$, yielding
\begin{equation}
t_{n}=\frac{1}{16}\left(\left(8+5\sqrt{2}\right)\left(1+\sqrt{2}\right)^{2n}+\left(8-5\sqrt{2}\right)\left(1-\sqrt{2}\right)^{2n}\right)-\frac{1}{2}
\end{equation}
Similarly, the closed forms for the three other variables are:

- for $\xi_{n}$ (\ref{eq:3.3-1}) (see OEIS \cite{key-13}, A336625):
\begin{equation}
\xi_{n}=\frac{1}{4}\left(\left(3-(-1)^{n}\sqrt{2}\right)\left(1+\sqrt{2}\right)^{2n+1}+\left(3+(-1)^{n}\sqrt{2}\right)\left(1-\sqrt{2}\right)^{2n+1}\right)-\frac{1}{2}
\end{equation}
yielding for even $n$ (upper sign) and odd $n$ (lower sign):
\begin{equation}
\xi_{n}=\frac{1}{4}\left(\left(3\mp\sqrt{2}\right)\left(1+\sqrt{2}\right)^{2n+1}+\left(3\pm\sqrt{2}\right)\left(1-\sqrt{2}\right)^{2n+1}\right)-\frac{1}{2}
\end{equation}
- for $T_{t_{n}}$ (\ref{eq:3.3-2}) (see OEIS \cite{key-13}, A336624):
\begin{align}
T_{t_{n}} & =\frac{1}{256}\left(\left(11-\left(-1\right)^{n}6\sqrt{2}\right)\left(1+\sqrt{2}\right)^{2\left(2n+1\right)}\right.\\
 & \left.+\left(11+\left(-1\right)^{n}6\sqrt{2}\right)\left(1-\sqrt{2}\right)^{2\left(2n+1\right)}\right)-\frac{9}{128}\nonumber 
\end{align}
yielding for even $n$ (upper sign) and odd $n$ (lower sign):
\begin{equation}
T_{t_{n}}=\frac{1}{256}\left(\left(3\mp\sqrt{2}\right)^{2}\left(1+\sqrt{2}\right)^{2\left(2n+1\right)}+\left(3\pm\sqrt{2}\right)^{2}\left(1-\sqrt{2}\right)^{2\left(2n+1\right)}\right)-\frac{9}{128}
\end{equation}
- for $T_{\xi_{n}}$ (\ref{eq:3.3-3}) (see OEIS \cite{key-13}, A336626):
\begin{align}
T_{\xi_{n}} & =\frac{1}{32}\left(\left(11\left(1+\sqrt{2}\right)^{2}-(-1)^{n}6\left(4+3\sqrt{2}\right)\right)\left(1+\sqrt{2}\right)^{4n}\right.\\
 & \left.+\left(11\left(1-\sqrt{2}\right)^{2}-(-1)^{n}6\left(4-3\sqrt{2}\right)\right)\left(1-\sqrt{2}\right)^{4n}\right)-\frac{9}{16}\nonumber 
\end{align}
yielding for even $n$:
\begin{equation}
T_{\xi_{n}}=\frac{1}{32}\left(\left(1+2\sqrt{2}\right)^{2}\left(1+\sqrt{2}\right)^{4n}+\left(1-2\sqrt{2}\right)^{2}\left(1-\sqrt{2}\right)^{4n}\right)-\frac{9}{16}
\end{equation}
and for odd $n$:
\begin{equation}
T_{\xi_{n}}=\frac{1}{32}\left(\left(5+4\sqrt{2}\right)\left(1+\sqrt{2}\right)^{4n}+\left(5-4\sqrt{2}\right)\left(1-\sqrt{2}\right)^{4n}\right)-\frac{9}{16}
\end{equation}

\subsection{Closed form equations for $r=3$ ($\rho=1$)}

The homogeneous solution (\ref{eq:5-10}) yields a complete solution
that reduces to 
\begin{align}
t_{n} & =\left(A+A^{\prime}\cos\left(\frac{2n\pi}{3}\right)+A^{\prime\prime}\sin\left(\frac{2n\pi}{3}\right)\right)\alpha^{n}\\
 & +\left(B+B^{\prime}\cos\left(\frac{2n\pi}{3}\right)+B^{\prime\prime}\sin\left(\frac{2n\pi}{3}\right)\right)\beta^{n}-1/2\nonumber \\
 & =\left(A\alpha^{n}+B\beta^{n}\right)+\left(A^{\prime}\alpha^{n}+B^{\prime}\beta^{n}\right)\cos\left(\frac{2n\pi}{3}\right)\nonumber \\
 & +\left(A^{\prime\prime}\alpha^{n}+B^{\prime\prime}\beta^{n}\right)\sin\left(\frac{2n\pi}{3}\right)-1/2\nonumber 
\end{align}
with the boundary conditions $t_{0}=t_{-1}=0$, $t_{1}=t_{-2}$ and
$t_{2}=t_{-3}$, giving six relations from which one finds the following
expressions of the six constants
\begin{align*}
A & =\frac{\alpha\left(\alpha^{2}-\beta^{3}\right)+\left(2t_{1}+1\right)\alpha^{2}\left(1-\beta^{3}\right)+\left(2t_{2}+1\right)\left(\alpha-1\right)}{6\left(\alpha^{3}-\beta^{3}\right)}\\
A^{\prime} & =\frac{\alpha\left(2\alpha+\beta^{3}\right)-\left(2t_{1}+1\right)\alpha^{2}\left(1-\beta^{3}\right)-\left(2t_{2}+1\right)\left(\alpha+2\right)}{6\left(\alpha^{3}-\beta^{3}\right)}\\
A^{\prime\prime} & =\frac{\sqrt{3}\alpha\left(\beta^{3}+\left(2t_{1}+1\right)\alpha\left(1-\beta^{3}\right)-\left(2t_{2}+1\right)\right)}{6\left(\alpha^{3}-\beta^{3}\right)}\\
B & =\frac{\beta\left(\alpha^{3}-\beta^{2}\right)+\left(2t_{1}+1\right)\beta^{2}\left(\alpha^{3}-1\right)-\left(2t_{2}+1\right)\left(\beta-1\right)}{6\left(\alpha^{3}-\beta^{3}\right)}\\
B^{\prime} & =\frac{-\beta\left(\alpha^{3}+2\beta^{2}\right)-\left(2t_{1}+1\right)\beta^{2}\left(\alpha^{3}-1\right)+\left(2t_{2}+1\right)\left(\beta+2\right)}{6\left(\alpha^{3}-\beta^{3}\right)}\\
B^{\prime\prime} & =\frac{\sqrt{3}\beta\left(-\alpha^{3}+\left(2t_{1}+1\right)\beta\left(\alpha^{3}-1\right)+\left(2t_{2}+1\right)\right)}{6\left(\alpha^{3}-\beta^{3}\right)}
\end{align*}
The smallest case of $k$ with rank three is $k=10$, for which $\kappa=18$
yielding $\alpha=\left(19+6\sqrt{10}\right)^{\frac{1}{3}}=\left(3+\sqrt{10}\right)^{\frac{2}{3}}$
and $\beta=\left(19-6\sqrt{10}\right)^{\frac{1}{3}}=\left(3-\sqrt{10}\right)^{\frac{2}{3}}$,
with $t_{1}=1$, $t_{2}=6$.

\subsection{Closed form equations for $r=4$ (\textmd{$\rho=2$)}}

The homogeneous solution (\ref{eq:5-11}) yields a complete solution
that reduces to 
\begin{align}
t_{n} & =\left(A+\left(-1\right)^{n}A^{\prime}+A^{\prime\prime}\cos\left(\frac{n\pi}{2}\right)+A^{\prime\prime\prime}\sin\left(\frac{n\pi}{2}\right)\right)\alpha^{n}\\
 & +\left(B+\left(-1\right)^{n}B^{\prime}+B^{\prime\prime}\cos\left(\frac{n\pi}{2}\right)+B^{\prime\prime\prime}\sin\left(\frac{n\pi}{2}\right)\right)\beta^{n}-1/2\nonumber 
\end{align}
and that simplifies as follows:

- for $n\equiv0\,\text{mod\ensuremath{4}}$, $t_{n}=\left(A+A^{\prime}+A^{\prime\prime}\right)\alpha^{n}+\left(B+B^{\prime}+B^{\prime\prime}\right)\beta^{n}-1/2$

- for $n\equiv1\,\text{mod\ensuremath{4}}$, $t_{n}=\left(A-A^{\prime}+A^{\prime\prime\prime}\right)\alpha^{n}+\left(B-B^{\prime}+B^{\prime\prime\prime}\right)\beta^{n}-1/2$

- for $n\equiv2\,\text{mod\ensuremath{4}}$, $t_{n}=\left(A+A^{\prime}-A^{\prime\prime}\right)\alpha^{n}+\left(B+B^{\prime}-B^{\prime\prime}\right)\beta^{n}-1/2$

- for $n\equiv3\,\text{mod\ensuremath{4}}$, $t_{n}=\left(A-A^{\prime}-A^{\prime\prime\prime}\right)\alpha^{n}+\left(B-B^{\prime}-B^{\prime\prime\prime}\right)\beta^{n}-1/2$

with the boundary conditions $t_{0}=t_{-1}=0$, $t_{1}=t_{-2}$, $t_{2}=t_{-3}$
and $t_{3}=t_{-4}$, giving eight relations yielding the eight constants

\begin{align*}
A & =\frac{\alpha\left(\alpha^{3}-\beta^{4}\right)+\left(2t_{1}+1\right)\alpha^{2}\left(\alpha-\beta^{4}\right)+\left(2t_{2}+1\right)\alpha^{2}\left(1-\alpha\beta^{4}\right)+\left(2t_{3}+1\right)\left(\alpha-1\right)}{8\left(\alpha^{4}-\beta^{4}\right)}\\
A^{\prime} & =\frac{\alpha\left(\alpha^{3}+\beta^{4}\right)-\left(2t_{1}+1\right)\alpha^{2}\left(\alpha+\beta^{4}\right)+\left(2t_{2}+1\right)\alpha^{2}\left(1+\alpha\beta^{4}\right)-\left(2t_{3}+1\right)\left(\alpha+1\right)}{8\left(\alpha^{4}-\beta^{4}\right)}\\
A^{\prime\prime} & =\frac{\alpha^{4}+\left(2t_{1}+1\right)\alpha^{2}\beta^{4}-\left(2t_{2}+1\right)\alpha^{2}-\left(2t_{3}+1\right)}{4\left(\alpha^{4}-\beta^{4}\right)}\\
A^{\prime\prime\prime} & =\frac{\alpha\left(\beta^{4}+\left(2t_{1}+1\right)\alpha^{2}-\left(2t_{2}+1\right)\alpha^{2}\beta^{4}-\left(2t_{3}+1\right)\right)}{4\left(\alpha^{4}-\beta^{4}\right)}\\
B & =\frac{\beta\left(\alpha^{4}-\beta^{3}\right)+\left(2t_{1}+1\right)\beta^{2}\left(\alpha^{4}-\beta\right)+\left(2t_{2}+1\right)\beta^{2}\left(\alpha^{4}\beta-1\right)-\left(2t_{3}+1\right)\left(\beta-1\right)}{8\left(\alpha^{4}-\beta^{4}\right)}\\
B^{\prime} & =\frac{-\beta\left(\alpha^{4}+\beta^{3}\right)+\left(2t_{1}+1\right)\beta^{2}\left(\alpha^{4}+\beta\right)-\left(2t_{2}+1\right)\beta^{2}\left(\alpha^{4}\beta+1\right)+\left(2t_{3}+1\right)\left(\beta+1\right)}{8\left(\alpha^{4}-\beta^{4}\right)}\\
B^{\prime\prime} & =\frac{-\beta^{4}-\left(2t_{1}+1\right)\alpha^{4}\beta^{2}+\left(2t_{2}+1\right)\beta^{2}+\left(2t_{3}+1\right)}{4\left(\alpha^{4}-\beta^{4}\right)}\\
B^{\prime\prime\prime} & =\frac{\beta\left(-\alpha^{4}-\left(2t_{1}+1\right)\beta^{2}+\left(2t_{2}+1\right)\alpha^{4}\beta^{2}+\left(2t_{3}+1\right)\right)}{4\left(\alpha^{4}-\beta^{4}\right)}
\end{align*}
The smallest case of $k$ with rank four is $k=13$, for which $\kappa=648$
yielding $\alpha=\left(649+180\sqrt{13}\right)^{\frac{1}{4}}=\left(18+5\sqrt{13}\right)^{\frac{1}{2}}$
and $\beta=\left(649-180\sqrt{13}\right)^{\frac{1}{4}}=\left(18-5\sqrt{13}\right)^{\frac{1}{2}}$,
with $t_{1}=3$, $t_{2}=21$ and $t_{3}=234$.

\section{Conclusions}

We have shown that for triangular numbers multiple of other triangular
numbers, closed form equations can be developed instead of recurrent
relations, to calculate directly the values of triangular numbers
and their indices without knowing the whole preceding sequence of
values. We developed the theoretical equations for four cases of ranks
from 1 to 4 and we applied them in several examples for the non-square
multipliers 2, 3, 5 and 8.

\end{document}